\documentclass{amsart}

\usepackage{amssymb}
\usepackage{amsmath}
\usepackage{amsthm}
\usepackage[english]{babel}
\usepackage[latin1]{inputenc}
\usepackage{mathrsfs}
\usepackage[all]{xy}

\newcommand{\F}{\ensuremath{\mathbf{F}}}

\newcommand{\h}[1]{\-\mbox{-#1}}

\newcommand{\I}{\ensuremath{\mathbf{I}}}

\newcommand{\kk}{\ensuremath{K}}

\newcommand{\OO}{\ensuremath{\mathcal{O}}}

\newcommand{\pp}{\ensuremath{\mathbb{P}}}

\newcommand{\syz}{\ensuremath{\mathrm{Syz}}}

\newcommand{\zz}{\ensuremath{\mathbb{Z}}}

\newtheorem{te}{Theorem}[section]
\newtheorem*{te*}{Theorem}
\newtheorem*{tecat*}{Teorema}
\newtheorem{p}[te]{Proposition}
\newtheorem*{p*}{Proposition}
\newtheorem*{pcat*}{Proposició}

\newtheorem*{co*}{Corollary}
\newtheorem*{cocat*}{Coro{\lgem}ari}

\newtheorem{lem}[te]{Lemma}
\newtheorem*{lem*}{Lemma}

\theoremstyle{definition}
\newtheorem{defin}[te]{Definition}

\newtheorem{pr}[te]{Problem}
\newtheorem*{pr*}{Problem}
\newtheorem*{prcat*}{Problema}
\newtheorem{st}[te]{Strategy}
\newtheorem{ob}[te]{Remark}

\title{Stability of syzygy bundles}

\date{\today}

\author{Pedro Macias Marques}
\address{Departamento de Matem\'atica, Universidade de \'Evora, Rua Rom\~{a}o Ramalho, 59, 7000--671 \'Evora, Portugal}
\curraddr{Departament d'\`Algebra i Geometria, Facultat de Matem\`atiques, Universitat de Barcelona, Gran Via de les Corts Catalanes, 585, 08007 Barcelona, Espanya}
\email{pmm@uevora.pt}
\urladdr{home.uevora.pt/$\sim$pmm}

\author{Rosa Mar\'\i a Mir\'o-Roig}
\address{Departament d'\`Algebra i Geometria, Facultat de Matem\`atiques, Universitat de\linebreak[4] Barcelona, Gran Via de les Corts Catalanes, 585, 08007 Barcelona, Espanya}
\email{miro@ub.edu}
\urladdr{atlas.mat.ub.es/personals/miro}

\thanks{The first author was partially supported by Funda\c{c}\~ao para a Ci\^encia e a Tecnologia, under grant SFRH/BD/27929/2006, and by CIMA -- Centro de Investiga\c{c}\~ao  em Matem\'atica e Aplica\c{c}\~oes, Universidade de \'Evora. The second author was partially supported by MTM2007-61104.}
\thanks{The authors wish to thank Holger Brenner for sharing his notes on the case of monomials in three variables. The idea behind lemma~\ref{Brenner} is his. The authors would also like to thank the referee, for the thorough proofreading and useful suggestions.}



\begin{document}

\begin{abstract}
We show that given integers~$N$, $d$ and~$n$ such that ${N\ge2}$,\linebreak[4] ${(N,d,n)\ne(2,2,5)}$, and ${N+1\le n\le\tbinom{d+N}{N}}$, there is a family of~$n$ monomials in $K\left[X_0,\ldots,X_N\right]$ of degree~$d$ such that their syzygy bundle is stable. Case ${N\ge3}$ was obtained independently by Coand\v{a} with a different choice of families of monomials \cite{Coa09}.

For ${(N,d,n)=(2,2,5)}$, there are~$5$ monomials of degree~$2$ in $K\left[X_0,X_1,X_2\right]$ such that their syzygy bundle is semistable.
\end{abstract}

\commby{Bernd Ulrich}

\maketitle

\section{Introduction}

Let $\kk$ be an algebraically closed field, ${R:=\kk[X_0,\ldots,X_N]}$ and ${\mathfrak{m}:=(X_0,\ldots,X_N)}$. A syzygy bundle is defined as the kernel of an epimorphism
\[\xymatrix{{\bigoplus\limits_{i=1}^n\OO_{\pp^N}(-d_i)}
    \ar[rr]^-{f_1,\ldots,f_n}&&\OO_{\pp^N}\mbox{,}}\]
given by ${(g_1,\ldots,g_n)\mapsto f_1g_1+\cdots+f_ng_n}$, where ${f_1,\ldots,f_n}$ are homogeneous polynomials in~$R$ of degrees ${d_1,\ldots,d_n}$, respectively, such that the ideal ${(f_1,\ldots,f_n)}$ is \mbox{$\mathfrak{m}$-primary}.

The main goal of this work is to give a complete answer to the following problem, presented by Brenner in \cite{Bre08b}:
\begin{pr}\label{6.9}
Does there exist for every~$d$ and every $n\le\tbinom{d+N}{N}$ a family of~$n$ monomials in~$R$ of degree~$d$ such that their syzygy bundle is semistable?
\end{pr}

In \cite{Bre08a}, Brenner computes the maximal slope of a syzygy bundle given by momnomials (see theorem~6.3). As a corollary, he deduces the following result (corollary~6.4, in his paper), which will be used as a main tool here:
\begin{p}\label{cor6.6}
Let $\left\{f_i\right\}_{i\in I}$ be a family of monomials in $R$ of degrees~$d_i$, such that the ideal ${(f_i,i\in I)}$ is \mbox{$\mathfrak{m}$-primary}. Suppose that, for every subset ${J\subseteq I}$, with ${|J|\ge2}$, the inequality
\[\frac{d_J-\sum_{i\in J}d_i}{|J|-1}\le\frac{-\sum_{i\in I}d_i}{|I|-1}\]
holds, where $d_J$ is the degree of the highest common factor of the subfamily $\left\{f_i\right\}_{i\in J}$. Then the syzygy bundle ${\syz(f_i,i\in I)}$ is semistable (and stable if strict inequality holds for ${J\subset I}$).
\end{p}

Note that when monomials are of the same degree~$d$, and making ${n:=|I|}$ and ${k:=|J|}$, the inequality in proposition~\ref{cor6.6} becomes
\begin{equation}\label{6.6}
    \left(d-d_J\right)n+d_J-dk\ge0\mbox{.}
\end{equation}

\smallskip

Case ${N=2}$ was solved in \cite{CMMR10}. We refer to this paper's introduction for more information on the problem and matters involved. The results in the present work are part of the first's author PhD thesis \cite{Mar09}. Case ${N\ge3}$ of the main result, theorem~\ref{main3}, was obtained independently by Coand\v{a} with a different choice of families of monomials \cite{Coa09}.

Let us now briefly explain how the paper is organized and how the families are constructed.

In section~\ref{N=1}, case ${N=1}$ is solved both applying direct methods and making use of the numerical criterion in proposition~\ref{cor6.6}. In section~\ref{N>2tools}, the tools to solve the general case ${N\ge3}$ are presented, and in section~\ref{N>2} its different subcases are dealt with. In the general case, we take the aproach described in the next paragraphs, based on what was done in~\cite{CMMR10}.

In general monomials in $R$ of a given degree~$d$ can be represented in a hypertetrahedron. This hypertetrahedron is the graph whose vertexes are all monomials of degree~$d$, and where two monomials are connected by an edge if and only if their greatest common divisor has degree ${d-1}$. We shall call the \emph{$i$th face} of this hypertetrahedron the set of monomials where the variable $X_i$ does not occur.

We will distinguish four cases, according to different values of~$n$. Recall that we have ${N+1\le n\le\tbinom{d+N}{N}}$. For the first cases given by
\[{N+1\le n\le\tbinom{d+N-1}{N-1}+1}\mbox{,}\]
we will show in lemma \ref{part of a face and opposite vertex} that each family of ${n-1}$ monomials in $\kk[X_0,\ldots,X_{N-1}]$ whose syzygy bundle over $\pp^{N-1}$ is stable yields a family of~$n$ monomials in $R$ whose syzygy bundle over $\pp^N$ is also stable. Cases
\[{\tbinom{d+N-1}{N-1}+1<n\le\tbinom{d+N}{N}-\tbinom{d-1}{N}}\]
are solved in proposition~\ref{faces} by taking the $N$th face and the vertex $X_N^d$ of the hypertetrahedron, and adding monomials in the remaining faces. Taking the set of all the hypertetrahedron's faces and adding the monomials in its interior which are closest to the vertexes gives us a solution to the cases
\[{\tbinom{d+N}{N}-\tbinom{d-1}{N}<n\le\tbinom{d+N}{N}
    -\tbinom{d-1}{N}+N+1}\mbox{,}\]
treated in proposition~\ref{facesanddots}. The last cases, with \[{\tbinom{d+N}{N}-\tbinom{d-1}{N}+N+1<n\le\tbinom{d+N}{N}}\mbox{,}\]
in lemma~\ref{Brenner}, are solved by taking a family of monomials of degree ${d-N-1}$ whose syzygy bundle is stable, multiplying them by ${X_0\cdots X_N}$, and adding all monomials in the faces of the hypertetrahedron. This is a generalisation of a lemma by Brenner, made for the case ${N=2}$ in his notes \cite{Bre}, which he kindly shared.

\section{Stable syzygy bundles on the projective line}
\label{N=1}

In this section a solution to problem~\ref{6.9} is presented for ${N=1}$. If ${n=2}$, we get stability, since the syzygy bundle is a line bundle. For ${n\ge3}$,  all vector bundles are a sum of line bundles, and therefore cannot be stable. We get thus the result:

\begin{te}\label{main1}
Let~$d$ and~$n$ be integers such that $2\le n\le d+1$ and~$d$ is a multiple of ${n-1}$. Then there is a family of~$n$ monomials in $K\left[X_0,X_1\right]$ of degree~$d$ such that their syzygy bundle is semistable. It is stable for ${n=2}$, and semistable, but not stable, otherwise. If~$d$ is not a multiple of ${n-1}$, there is no such family.

Moreover, if~$d$ is not a multiple of ${n-1}$, and ${f_1,\ldots,f_n}$ is any family of homogeneous polynomials in $\kk[X_0,X_1]$ of degree~$d$ such that the ideal ${(f_1,\ldots,f_n)}$ is \mbox{$\mathfrak{m}$-primary}, their syzygy bundle ${\syz(f_1,\ldots,f_n)}$ is not semistable.
\end{te}
\begin{proof}
If ${n=d+1}$, if  $\I:=\big\{X_0^d,\,X_0^{d-1}X_1,\ldots,X_1^d\big\}$ then the syzygy bundle $\syz(\I)$ is semistable, but not stable. Indeed, if~$g$ is the greatest common divisor of monomials in a subset ${J\subseteq\I}$, all monomials in~$J$ are of the form $gh$, with~$h$ a monomial of degree ${d-d_J}$, where~$d_J$ is the degree of~$g$. There are ${d-d_J+1}$ monomials of degree ${d-d_J}$, so ${k:=|J|\le d-d_J+1}$. Now
\[(d-d_J)n+d_J-dk\ge(d-d_J)(d+1)+d_J-d(d-d_J+1)=0.\]
Therefore inequality~(\ref{6.6}) holds. In fact, if we consider the subfamily of monomials
$J:=\big\{X_0^d,\,X_0^{d-1}X_1,
    \ldots,X_0^{d-d_J}X_1^{d_J}\big\}$,
we get equality, which means that the syzygy bundle is not stable.

In general, for the remaining values of~$n$, i.e.\ $3\le n\le d$, if ${f_1,\ldots,f_n}$ is a family of homogeneous polynomials in $\kk[X_0,X_1]$ such that the ideal ${(f_1,\ldots,f_n)}$ is \mbox{$\mathfrak{m}$-primary}, their syzygy bundle ${\syz(f_1,\ldots,f_n)}$ has rank ${n-1}$ and first Chern class ${c_1\big(\syz(f_1,\ldots,f_n)\big)=-dn}$. By Grothendieck theorem, there are integers ${a_1,\ldots,a_{n-1}}$ such that
${\syz(f_1,\ldots,f_n)\cong
    \bigoplus_{i=1}^{n-1}\OO_{\pp^N}(a_i)}$
which is semistable if and only if $a_1=\cdots=a_{n-1}$. Therefore
\[(n-1)a_1=a_1+\cdots+a_{n-1}=
    c_1\big(\syz(f_1,\ldots,f_n)\big)=-dn\mbox{,}\]
and since $n$ and ${n-1}$ are coprime, $d$ is a multiple of ${n-1}$. We have thus found a necessary condition for such a syzygy bundle to be semistable.

Now for the converse, suppose $d$ is a multiple of ${n-1}$, say ${d=(n-1)e}$, with ${e\in\zz}$, and consider the set
\[\I:=\big\{X_0^d,\,X_0^{(n-2)e}X_1^e,
    \,X_0^{(n-3)e}X_1^{2e},\ldots,\,X_0^eX_1^{(n-2)e},
    \,X_1^d\big\}.\]
Here we can get an isomorphism $\syz(\I)\cong\OO_{\pp^N}(-ne)^{n-1}$ by sending standard vectors to $\big(0,\ldots,X_1^e,\,-X_0^e,\ldots,0\big)$. Therefore, $\syz(\I)$ is a semistable bundle. Since in cases ${n=2}$ and ${n=d+1}$, $d$ is also a multiple of ${n-1}$, we have the result.
\end{proof}

\section{Some tools for the general case}
\label{N>2tools}

As we have said in the introduction, case ${N=2}$ was solved in \cite{CMMR10}. Let us now consider case ${N\ge3}$. We start with some tools and definitions we will need, and deal with the several cases in the next section.
\begin{ob}\label{adj}
If we use the notation $a_{d,j}:=-\tfrac{jd}{j-1}$, inequality~(\ref{6.6}) is equivalent to ${\tfrac{d_J}{k-1}+a_{d,k}\le a_{d,n}}$. The fact that once~$d$ is fixed, the sequence $\left(a_{d,j}\right)_{j\ge2}$ is monotonically increasing will be useful in many arguments.
\end{ob}
\begin{lem}\label{part of a face and opposite vertex}
If $N\ge3$, $N+1\le n\le\tbinom{d+N-1}{N-1}+1$, and $\I'$ is a family of ${n-1}$ monomials in $\kk\left[X_0,\ldots,X_{N-1}\right]$ of degree~$d$ such that their syzygy bundle is stable, then  ${\I:=\I'\cup\big\{X_N^d\big\}}$ is a family of~$n$ monomials in $R$ of degree~$d$ whose associated syzygy bundle is stable.
\end{lem}
\begin{proof}
Note that the ideal generated by~$\I$ is primary. Let ${J\subseteq\I}$ be a subset with at least two monomials. If ${J\subseteq \I'}$, then by hypotheses, inequality~(\ref{6.6}) holds. If not, then $X_N^d\in J$, and since~$J$ has at least another monomial, where the variable~$X_N$ does not occur and $d_J=0$, so inequality~(\ref{6.6}) holds, for the sequence $\left(a_{d,j}\right)_{j\ge2}$ is monotonically increasing.
\end{proof}

A direct application of this lemma which will become handy for proofs to follow is to take the whole $N$th face and add $X_N^d$ to obtain a well\h{behaved} family again. To this end, let us state the case of the highest possible~$n$, which has already been proved by Flenner in characteristic zero \cite{Fle84}, by Ballico \cite{Bal92}, and by Brenner \cite{Bre08a}, using his own criterion.

\begin{p}\label{casemaximalngeneral}
For any ${N\ge2}$,
$\syz\left(\big\{X_0^{i_0}\cdots X_N^{i_N}:
    i_0+\cdots+i_N=d\big\}\right)$
is stable on~$\pp^N$.
\end{p}
\begin{proof}
Let $\I:=\big\{X_0^{i_0}\cdots X_N^{i_N}: i_0+\cdots+i_N=d\big\}$. If~$g$ is the greatest common divisor of monomials in a subset ${J\subseteq \I}$, all monomials in~$J$ are of the form $gh$, with~$h$ a monomial of degree ${d-d_J}$, where~$d_J$ is the degree of~$g$. There are $\tbinom{N+d-d_J}{d-d_J}$ monomials of degree ${d-d_J}$, so ${k=|J|\le\tbinom{N+d-d_J}{d-d_J}}$. Now
\[(d-d_J)n+d_J-dk=(d-d_J)\tbinom{d+N}{N}+d_J
    -d\tbinom{d-d_J+N}{N}>0,\]
which can be proved by induction on~$N$. Therefore inequality~(\ref{6.6}) holds.
\end{proof}

Now using lemma~\ref{part of a face and opposite vertex} and this
proposition, we get

\begin{p}\label{whole face and opposite vertex}
For any ${N\ge3}$, the syzygy bundle associated to the family
$\big\{X_0^{i_0}\cdots X_{N-1}^{i_{N-1}}:
    i_0+\cdots+i_{N-1}=d\big\}\cup\big\{X_N^d\big\}$
is stable on~$\pp^N$.
\end{p}

In what follows, except when stated otherwise, we will adopt the
following:

\begin{st}\label{X_0}
For each given~$d$ and~$n$, we choose a set of~$n$ monomials~$\I$ such that for ${0<d_J<d}$, no monomial of degree~$d_J$ divides a greater number of monomials in~$\I$ than~$X_0^{d_J}$. 
\end{st}
\begin{ob}\label{maxJ}
If ${J\subset\I}$, with ${k:=|J|\ge2}$, and $d_J$ is the degree of the greatest common divisor of monomials in~$J$, to verify that~$\I$ satisfies inequality~(\ref{6.6}), we can assume ${0<d_J<d}$, since the fact that~$J$ has at least two elements makes ${d_J\ne d}$, and for ${d_J=0}$, the fact that the sequence $\left(a_{d,j}\right)_{j\ge2}$ is monotonically increasing is enough. We may also assume that~$J$ has all multiples of its greatest common divisor, since if a degree~$d_J$ is fixed, the higher~$k$ is, the harder it is to guarantee inequality~(\ref{6.6}).
\end{ob}

Finally let us define the set of all faces of the hypertetrahedron. 

\begin{defin}\label{deffaces}
Let $\F_{N,d}$ denote the family of monomials
\[\F_{N,d}:=\big\{X_0^{i_0}\cdots X_N^{i_N}:
    i_0+\cdots+i_N=d\mbox{ and }i_0\cdots i_N=0\big\}\mbox{.}\]
\end{defin}
We can easily see that the cardinality of~$\F_{N,d}$ is $\tbinom{d+N}{N}-\tbinom{d-1}{N}.$

\section{Stable syzygy bundles on the projective space}
\label{N>2}

In this section, we use the tools of the previous one to prove our main result (theorem~\ref{main3}).

\begin{p}\label{faces}
Let ${N\ge3}$, ${d\ge2}$ and ${\tbinom{d+N-1}{N-1}+1<n\le\tbinom{d+N}{N}-\tbinom{d-1}{N}}$. Then there is a family of~$n$ monomials in $R$ of degree~$d$ whose associated syzygy bundle is stable.
\end{p}
\begin{proof}
The upper bound considered for $n$ in this proposition becomes $\tbinom{d+N}{N}$ when ${d-1<N}$. In case ${n=\tbinom{d+N}{N}}$ we know the statement is true by proposition~\ref{casemaximalngeneral}. Otherwise, let ${1\leq r\leq\min(d-1,N)}$ and ${0\leq l\leq d-r-1}$ be such that
\[\tbinom{d+N}{N}-\tbinom{d-r+N}{N}+\tbinom{l+N-1}{N-1}<n
    \leq\tbinom{d+N}{N}-\tbinom{d-r+N}{N}+\tbinom{l+N}{N-1}\mbox{,}\]
and let
$I'_r:=\big\{X_0^{j_0}\cdots X_N^{j_N}:
        j_0+\cdots+j_N=d\mbox{ and }
        j_{N-r+1}\cdots j_N=0\big\}$.
This set contains all monomials in faces  ${N-r+1}$ to~$N$, and $\left|I'_r\right|=\tbinom{d+N}{N}-\tbinom{d-r+N}{N}$. Now let
$I''_{r,l}$ be the set of the monomials in face ${N-r}$ with degree in $X_N$ greater than ${d-r-l}$ that do not belong to $I'_r$, i.e.\ the ones of type
\[{X_{N-r+1}\cdots X_{N-1}X_N^{d-r-l+1}f}\mbox{,}\]
where~$f$ is a monomial of degree ${l}$ where the variable $X_{N-r}$ does not occur. Therefore $\left|I''_{r,l}\right|=\tbinom{l+N-1}{N-1}$. Let ${1\le i
\le\tbinom{l+N-1}{N-2}}$ be such that
\[n=\tbinom{d+N}{N}-\tbinom{d-r+N}{N}+\tbinom{l+N-1}{N-1}+i\mbox{,}\]
and let~$I'''_{r,l}$ be a set of~$i$ monomials of degree~$d$ in $R$ of the form
\[{X_{N-r+1}\cdots X_{N-1}X_N^{d-r-l}f}\mbox{,}\]
where~$f$ is a monomial of degree~${l+1}$, where variables $X_{N-r}$ and $X_N$ do not occur. Let us choose these monomials in such a way that the degrees of~$X_0$ in these~$f$ are as large as possible. Let ${\I:=I'_r\cup I''_{r,l}\cup I'''_{r,l}}$. Since~${I'_r\cup I''_{r,l}}$ is a set for which the claim in strategy~\ref{X_0} is true, the way we choose the monomials for~$I'''_{r,l}$ guarantees that strategy~\ref{X_0} can be applied to~$\I$.

As always, it is enough to verify inequality~(\ref{6.6}) for ${0<d_J<d}$ (see strategy~\ref{X_0}). We shall see the cases ${0<d_J\le l}$, ${d_J=l+1}$ and ${l+1<d_J<d}$ separately.

\bigskip
\noindent{\sc Case 1:} ${d_J\leq l}$. In this case, if~$k$ is the
number of multiples of $X_0^{d_J}$ in~$\I$, we have
\[k=\tbinom{d-d_J+N}{N}-\tbinom{d-d_J-r+N}{N}
    +\tbinom{l-d_J+N-1}{N-1}+\min\left[i,\tbinom{l-d_J+N-1}{N-2}\right].\]
Therefore
\begin{multline*}
(d-d_J)n+d_J-dk\ge(d-d_J)\left[\tbinom{d+N}{N}
        -\tbinom{d-r+N}{N}+\tbinom{l+N-1}{N-1}\right]+d_J\\
    -d\left[\tbinom{d-d_J+N}{N}
        -\tbinom{d-d_J-r+N}{N}
        +\tbinom{l-d_J+N-1}{N-1}\right] -d_J\tbinom{l-d_J+N-1}{N-2}.
\end{multline*}
Concluding this case amounts to showing that this last expression is positive. To do this, set
\begin{align*}
T(N,d,d_J,r,l)&:=(d-d_J)\left[\tbinom{d+N}{N}
        -\tbinom{d-r+N}{N}+\tbinom{l+N-1}{N-1}\right]+d_J\\
    &\qquad -d\left[\tbinom{d-d_J+N}{N}
        -\tbinom{d-d_J-r+N}{N}
        +\tbinom{l-d_J+N-1}{N-1}\right]-d_J\tbinom{l-d_J+N-1}{N-2}.
\end{align*}

Let us start by showing that~$T$ increases with~$r$:
\[T(N,d,d_J,r+1,l)-T(N,d,d_J,r,l)=
    (d-d_J)\tbinom{d-r+N-1}{N-1}-d\tbinom{d-d_J-r+N-1}{N-1}\]
Note that if ${r>d-d_J}$, then
${\tbinom{d-d_J-r+N-1}{N-1}=0\mbox{,}}$ in which case this
expression is clearly positive. Otherwise, we get
{\small
\begin{align*}
&T(N,d,d_J,r+1,l)-T(N,d,d_J,r,l)=\\
    &=\tfrac{1}{(N-1)!}\left[(d-d_J)(d-r+1)
        \prod_{s=2}^{N-1}(d-r+s)
        -d(d-d_J-r+1)
        \prod_{s=2}^{N-1}(d-d_J-r+s)\right].
\end{align*}
}%
This last expression is non\h{negative}, since
\[(d-d_J)(d-r+1)-d(d-d_J-r+1)=(r-1)d_J\ge0\mbox{,}\]
and ${d-r+s>d-d_J-r+s}$. We can therefore look at the case ${r=1}$, since if~$T$ is positive in this case, it will always be positive. Let us see now that~$T(N,d,d_J,1,l)$ increases with~$l$. Suppose ${d_J\le l\le d-3}$. We get
{\small
\begin{multline*}
T(N,d,d_J,1,l+1)-T(N,d,d_J,1,l)=\\
    =\frac{1}{(N-2)!}\left[(d-d_J)\prod_{s=2}^{N-1}(l+s) -\big[d(l-d_J+2)+d_J(N-2)\big]
        \prod_{s=2}^{N-2}(l-d_J+s+1)\right].
\end{multline*}
}%
This last expression is never negative, since for ${N=3}$ we have
\[T(3,d,d_J,1,l+1)-T(3,d,d_J,1,l)=(d-l-3)d_J\ge0,\]
and for ${N\ge4}$ we can write
{\small
\begin{multline*}
T(N,d,d_J,1,l+1)-T(N,d,d_J,1,l)=
    \frac{1}{(N-2)!}\left[(d-d_J)(l+2)(l+N-1)
        \prod_{s=3}^{N-2}(l+s)\right. \\
    \left.-\big[d(l-d_J+2)+d_J(N-2)\big]
        (l-d_J+3)\prod_{s=3}^{N-2}(l-d_J+s+1)\right]
\end{multline*}
}%
This is not negative, since
\begin{align*}
(d-d_J)&(l+N-1)(l+2)-\big[d(l-d_J+2)+d_J(N-2)\big](l-d_J+3)=\\
    &=(d-l-3)l(N-2)+(l-d_J)^2(N-3)+2(d-l)(N-3)\\
    &\qquad +5(l-d_J)(N-3)+2(d-l-1)(l-d_J)(d_J-1)\\
    &\qquad +(d-l-3){d_J}^2+l(l-d_J)(d_J-1)+3(d-l-2)(d_J-1)\\
    &\qquad +(l-d_J)(d_J-1)+(d-l-3)+3(d_J-1)(d_J+1)\mbox{,}
\end{align*}
which is non\h{negative}, and for ${3\le s\le N-2}$ we have
$l+s\ge l-d_J+s+1.$

Therefore we can look at the case ${l=d_J}$, for if $T$ is positive in this case, it will always be positive. We look at two cases separately: ${d>2d_J}$ and ${d\le2d_J}$. In the former, we get
{\small
\begin{align*}
&T(N,d,d_J,1,d_J)=(d-d_J)\left[\tbinom{d+N-1}{N-1}
        +\tbinom{d_J+N-1}{N-1}\right]-d\tbinom{d-d_J+N-1}{N-1}-d-d_J(N-2)\\
    &=\tfrac{1}{(N-1)!}\left[(d-d_J)(d+1)(d+2)
        \prod_{s=3}^{N-1}(d+s)
   -d(d-d_J+1)(d-d_J+2)
        \prod_{s=3}^{N-1}(d-d_J+s)\right.\\
    &\qquad\left.{}+(d-2d_J)\left(
        \prod_{s=1}^{N-1}(d_J+s)-(N-1)!\right) +d_J\left(\prod_{s=1}^{N-1}(d_J+s)
        -N!\right)\right].
\end{align*}
}%
Now
\begin{multline*}
(d-d_J)(d+1)(d+2)-d(d-d_J+1)(d-d_J+2)=\\
    =(d-2d_J)d\cdot d_J+(d-2d_J){d_J}^2+2d_J\left({d_J}^2-1\right),
\end{multline*}
which is never negative, and we always have ${d+s>d-d_J+s}$. Furthermore
we can see that ${\prod_{s=1}^{N-1}(d_J+s)-N!\ge\prod_{s=1}^{N-1}(1+s)-N!=0}$, and the term\linebreak[4] ${(d-2d_J)\left[\prod_{s=1}^{N-1}(d_J+s)-(N-1)!\right]}$ is strictly positive. Therefore~$T$ is positive in this case.

In case ${d\le2d_J}$, since ${1\le d_J=l\le d-2}$, we get ${d\ge3}$, and therefore ${d_J\ge2}$. We can write
{\small
\begin{align*}
&T(N,d,d_J,1,d_J)=(d-d_J)\left[\tbinom{d+N-1}{N-1}
        +\tbinom{d_J+N-1}{N-1}\right]-d\tbinom{d-d_J+N-1}{N-1} -d-d_J(N-2)\\
    &=\tfrac{1}{(N-1)!}\left[(d-d_J)(d+1)(d+2)
        \prod_{s=3}^{N-1}(d+s)-d(d-d_J+1)(d-d_J+2)
        \prod_{s=3}^{N-1}(d-d_J+s)\right.\\
    &\qquad\qquad\qquad\left.{}+(d-d_J)(d_J+1)
        \prod_{s=3}^N(d_J-1+s)
        -2d_J\prod_{s=3}^Ns\right]+2d_J-d.
\end{align*}
}%
Now if we observe that
\begin{multline*}
(d-d_J)(d+1)(d+2)-d(d-d_J+1)(d-d_J+2)=\\
    =d(d-d_J-1)d_J+(d-2)d_J>0\mbox{,}
\end{multline*}
we see that the difference
\[(d-d_J)(d+1)(d+2)\prod_{s=3}^{N-1}(d+s)-d(d-d_J+1)(d-d_J+2)
    \prod_{s=3}^{N-1}(d-d_J+s)\]
is strictly positive. Finally, since ${d-d_J\ge2}$, we get $(d-d_J)(d_J+1)-2d_J\ge2(d_J+1)-2d_J>0\mbox{,}$ and
\[(d-d_J)(d_J+1)\prod_{s=3}^N(d_J-1+s)-2d_J\prod_{s=3}^Ns >0.\]

Therefore $T$ is positive also in this case, and hence it is always positive.

\bigskip
\noindent{\sc Case 2:} ${d_J=l+1}$. To count the multiples of $X_0^{d_J}$, i.e.\ $X_0^{l+1}$, in $I'$, we count all possible multiples of $X_0^{l+1}$ and subtract the ones of the form ${X_0^{l+1}fX_{N-r+1}\cdots X_N}$,
where~$f$ is a monomial of degree ${d-l-1-r}$. We get ${\tbinom{d-l-1+N}{N}-\tbinom{d-l-1-r+N}{N}}$. In $I''$ there is only one multiple of $X_0^{d_J}$, namely ${X_0^{l+1}X_{N-r+1}\cdots X_{N-1}X_N^{d-r-l}}$. Therefore if~$k$ is the number of multiples of $X_0^{d_J}$ in~$\I$, we have
\[k=\tbinom{d-l-1+N}{N}-\tbinom{d-l-1-r+N}{N}+1.\]
we get
\begin{align*}
(d-d_J)n+d_J-dk&=(d-l-1)\left[\tbinom{d+N}{N}
        -\tbinom{d-r+N}{N}+\tbinom{l+N-1}{N-1}+i\right]\\
    &\qquad +l+1-d\left[\tbinom{d-l-1+N}{N}
        -\tbinom{d-l-1-r+N}{N}+1\right]\displaybreak[0]\\
    &\ge(d-l-1)\left[\tbinom{d+N}{N}
        -\tbinom{d-r+N}{N}+\tbinom{l+N-1}{N-1}+1\right]\\
    &\qquad +l+1-d\left[\tbinom{d-l-1+N}{N}
        -\tbinom{d-l-1-r+N}{N}+1\right]\displaybreak[0]\\
    &=(d-l-1)\left[\tbinom{d+N}{N}
        -\tbinom{d-r+N}{N}+\tbinom{l+N-1}{N-1}\right]\\
    &\qquad -d\left[\tbinom{d-l-1+N}{N}
        -\tbinom{d-l-1-r+N}{N}\right]
\end{align*}
Again all we have to do is to prove that this expression is positive. Let
{\small
\[U(N,d,r,l):=(d-l-1)\left[\tbinom{d+N}{N}
        -\tbinom{d-r+N}{N}+\tbinom{l+N-1}{N-1}\right]-d\left[\tbinom{d-l-1+N}{N}
        -\tbinom{d-l-1-r+N}{N}\right],\]
}%
and, as we did in the previous case, let us see that this expression increases with~$r$.
{\small
\begin{align*}
U&(N,d,r+1,l)-U(N,d,r,l)=\\
    &=\tfrac{1}{(N-1)!}\left[(d-l-1)(d-r+1)
        \prod_{s=2}^{N-1}(d-r+s)-d(d-l-r)
        \prod_{s=2}^{N-1}(d-l-1-r+s)\right].
\end{align*}
}%
This is never negative, since
$(d-l-1)(d-r+1)-d(d-l-r)=(l+1)(r-1)\ge0$, and for ${2\le s\le
N-1}$, ${d-r+s\ge d-l-1-r+s}$. Therefore, to see that~$U$ is
positive, it is enough to check the case ${r=1}$.
\begin{multline*}
U(N,d,1,l)=(d-l-1)\left[\tbinom{d+N-1}{N-1}+\tbinom{l+N-1}{N-1}\right]
       -d\tbinom{d-l+N-2}{N-1}\\
    =\tfrac{1}{(N-1)!}\left[(d-l-1)(d+1)(d+2)
        \prod_{s=3}^{N-1}(d+s)\right]\\
     -\tfrac{1}{(N-1)!}\left[d(d-l)(d-l+1)
        \prod_{s=3}^{N-1}(d-l-1+s)-(d-l-1)\prod_{s=1}^{N-1}(l+s)\right]
\end{multline*}
Since
\begin{multline*}
(d-l-1)(d+1)(d+2)-d(d-l)(d-l+1)=\\
    =d(d-l-2)l+d(d-l-2)+(d-2)l+(d-2)\ge0\mbox{,}
\end{multline*}
and ${(d-l-1)\prod_{s=1}^{N-1}(l+s)>0}$, we get that~$U$ is positive.

\bigskip
\noindent{\sc Case 3:} ${d_J>l+1}$. In this case, we shall use induction on~$r$. For ${r=1}$, we get that all multiples of $X_0^{d_J}$ in~$\I$ are in~$I'_1$, and since~${I'_1\cup I''_{1,0}}$ is the set mentioned in proposition~\ref{whole face and opposite vertex}, inequality~(\ref{6.6}) is satisfied for~${I'_1\cup I''_{1,0}}$. The fact that the sequence $\left(a_{d,j}\right)_{j\ge2}$ is monotonically increasing guarantees that inequality~(\ref{6.6}) is also satisfied for~$\I$.

For the induction step, suppose that for a given~$r$, inequality~(\ref{6.6}) is satisfied for~${I'_r\cup I''_{r,d-r-1}\cup I'''_{r,d-r-1}}$. Note that ${I'_{r+1}\cup I''_{r+1,0}=I'_r\cup I''_{r,d-r-1}\cup I'''_{r,d-r-1}}$, if we consider the maximum possible~$i$, so that ${I'''_{r,d-r-1}}$ has cardinality ${\tbinom{d-r+N-2}{N-2}}$. Then again the fact that the sequence $\left(a_{d,j}\right)_{j\ge2}$ is monotonically increasing guarantees that inequality~(\ref{6.6}) is also satisfied for~$\I$ (see remark~\ref{adj}), since all multiples of $X_0^{d_J}$ belong to ${I'_{r+1}}$.
\end{proof}

When ${d=N+1}$, the proposition above leaves out only one monomial, namely ${X_0\cdots X_N}$. Proposition~\ref{casemaximalngeneral} shows that for ${n=\tbinom{d+N}{N}=\tbinom{2N+1}{N}}$, the whole hypertetrahedron is a family whose associated syzygy bundle is stable. For ${d>N+1}$, the next proposition starts with $\F_{N,d}$ (recall definition~\ref{deffaces}) and adds at most ${N+1}$ monomials.

\begin{p}\label{facesanddots}
Let ${N\ge3}$, ${d>N+1}$ and
\[\tbinom{d+N}{N}-\tbinom{d-1}{N}<n\le
    \tbinom{d+N}{N}-\tbinom{d-1}{N}+N+1.\]
Then there is a family of~$n$ monomials in $R$ of degree~$d$ whose associated syzygy bundle is stable.
\end{p}
\begin{proof}
Let ${1\le i\le N+1}$ be such that ${n=\tbinom{d+N}{N}-\tbinom{d-1}{N}+i}$,
and let~$I'$ be the set of the first~$i$ monomials in the sequence
\[\big(X_0^{d-N}X_1\cdots X_N,\,X_0X_1^{d-N}X_2\cdots X_N,
    \ldots,X_0\cdots X_{N-1}X_N^{d-N}\big).\]
Let ${\I:=\F_{N,d}\cup I'}$. Let us check that inequality~(\ref{6.6}) holds for ${0<d_J<d}$. Since~$\I$ satisfies strategy~\ref{X_0}, we can look only at multiples of~${X_0^{d_J}}$.

If ${d-N<d_J<d}$, all multiples of~${X_0^{d_J}}$ belong to~$\F_{N,d}$, and since, by the previous proposition, inequality~(\ref{6.6}) holds for~$\F_{N,d}$, the fact that the sequence $\left(a_{d,j}\right)_{j\ge2}$ is monotonically increasing guarantees that inequality~(\ref{6.6}) holds for~$\I$.

If $1<d_J\le d-N$ the number of multiples of~${X_0^{d_J}}$ is
$k:=\tbinom{d-d_J+N}{N}-\tbinom{d-d_J}{N}+1$, and we get
{\small
\begin{align*}
(&d-d_J)n+d_J-dk=(d-d_J)\left[\tbinom{d+N}{N}-\tbinom{d-1}{N}+i
        \right]+d_J-d\left[\tbinom{d-d_J+N}{N}
        -\tbinom{d-d_J}{N}+1\right]\\
    &\ge(d-d_J)\left[\tbinom{d+N}{N}-\tbinom{d-1}{N}+1\right]+d_J-d\left[\tbinom{d-d_J+N}{N}
        -\tbinom{d-d_J}{N}+1\right]\\
    &=(d-d_J)\left[\tbinom{d+N}{N}-\tbinom{d-1}{N}\right]
        -d\left[\tbinom{d-d_J+N}{N}
        -\tbinom{d-d_J}{N}\right].
\end{align*}
}%
Let
\[V(d,d_J,N):=(d-d_J)\left[\tbinom{d+N}{N}-\tbinom{d-1}{N}\right]
    -d\left[\tbinom{d-d_J+N}{N}-\tbinom{d-d_J}{N}\right].\]
If we look at~$V$ as a function on~$d_J$, its second derivative is
\[-\tfrac{2d}{N!}\cdot\sum_{1\le s<t\le N}\
    \left[\sideset{}{_{t\neq r\neq s}}\prod_{r=1}^N(d-d_J+r)
    -\sideset{}{_{t\neq r\neq s}}\prod_{r=1}^N
    (d-d_J-N+r)\right].\]
Since this is negative for ${d_J\in[1,d-N]}$, the function's minimum in this interval is at one of its ends. (Note that we are dealing with the case ${1<d_J\le d-N}$, and therefore the lowest value for $d_J$ is~$2$, and not~$1$, but for the sake of simplicity in calculations, we can look at ${1\le d_J\le d-N}$.) For ${d_J=1}$, we get
{\small
\[V(d,1,N)=\tfrac{1}{N!}\left[(d-1)(d+1)(d+2)\prod_{s=2}^N(d+s)
        -d^2(d+1)\prod_{s=2}^N(d-1+s)\right]+\tbinom{d-1}{N}.\]
}%
This is clearly positive, since
${(d-1)(d+1)(d+2)-d^2(d+1)=(d-2)(d+1)>0}$. For ${d_J=d-N}$, we get
{\small
\begin{multline*}
V(d,\,d-N,\,N)=\tfrac{1}{N!}\left[N(d+1)(d+2)
        \prod_{s=3}^N(d+s)\right]\\
    -\tfrac{1}{N!}\left[N(d-N)(d-N+1)\prod_{s=3}^N(d-N-1+s)
        +d(N+1)(N+2)\prod_{s=3}^N(N+s)\right]+d.
\end{multline*}
}%
This is easily seen to be positive, since
\begin{multline*}
N(d+1)(d+2)-N(d-N)(d-N+1)-d(N+1)(N+2)=\\
    =(d-N)N(N-2)+(d-N)(N-2)>0,
\end{multline*}
which means that ${N(d+1)(d+2)>N(d-N)(d-N+1)+d(N+1)(N+2)}$,
and for $3\le s\le N$, ${d+s>d-N-1+s}$ and ${d+s>N+s}$.

Finally, if $d_J=1$, the number of multiples of~${X_0^{d_J}}$ is
$k:=\tbinom{d-1+N}{N}-\tbinom{d-1}{N}+i$, and we get
\begin{align*}
(&d-d_J)n+d_J-dk=(d-1)\left[\tbinom{d+N}{N}-\tbinom{d-1}{N}+i
        \right] +1-d\left[\tbinom{d-1+N}{N}
        -\tbinom{d-1}{N}+i\right]\displaybreak[0]\\
    &=(d-1)\left[\tbinom{d+N}{N}-\tbinom{d-1}{N}\right]-d\left[\tbinom{d-1+N}{N}
        -\tbinom{d-1}{N}\right]+1-i\displaybreak[0]\\
    &\ge(d-1)\left[\tbinom{d+N}{N}-\tbinom{d-1}{N}\right]
        -d\left[\tbinom{d-1+N}{N}-\tbinom{d-1}{N}\right]-N\displaybreak[0]\\
    &=d\tbinom{d+N-1}{N-1}-\tbinom{d+N}{N}+\tbinom{d-1}{N}-N=\left(d-\tfrac{d+N}{N}\right)\tbinom{d+N-1}{N-1}
        +\tbinom{d-1}{N}-N.
\end{align*}
This is positive, since ${d-\tfrac{d+N}{N}=\tfrac{(d-2)(N-1)+N-2}{N}>0}$, and ${\tbinom{d-1}{N}\ge\tbinom{N+1}{N}>N}$.

Therefore inequality~(\ref{6.6}) strictly holds in all cases, so the syzygy bundle associated to~$\I$ is stable.
\end{proof}

For some computations in the next proof, we will need a result which is a simple consequence of the fact that for any numbers $a,b_1,\ldots,b_n$,
\[\prod_{s=1}^p(a+b_s)=\sum_{s=0}^p\
    \sum_{1\le t_1<\cdots<t_s\le p}
    a^{p-s}b_{t_1}\cdots b_{t_1}.\]
Note that if ${N\ge1}$ and ${d\ge0}$, we get ${\tbinom{d+N}{N}-\tbinom{d-1}{N}=\tfrac{N+1}{(N-1)!}d^{N-1}+
    \mbox{positive terms}}$.
\begin{lem}\label{Brenner2}
Let ${N\ge1}$ and ${d\ge0}$. Then
\[\tbinom{d+N}{N}-\tbinom{d-1}{N}\ge\tfrac{N+1}{(N-1)!}d^{N-1}.\]
\end{lem}

The following lemma will be the key to prove our main result.

\begin{lem}\label{Brenner}
Let ${N\ge3}$, ${d>N+1}$, and ${\tbinom{d+N}{N}-\tbinom{d-1}{N}<n\le\tbinom{d+N}{N}}$. If~$I'$ is a family of
\[{n':=n-\left[\tbinom{d+N}{N}-\tbinom{d-1}{N}\right]}\]
monomials in $R$ of degree~${d':=d-N-1}$ such that their syzygy bundle is stable, then
\[\I:=\F_{N,d}\cup\big\{X_0\cdots X_Nf:f\in I'\big\}\]
is a family of $n$ monomials in $R$ of degree~$d$ whose associated syzygy bundle is stable.
\end{lem}
\begin{proof}
Let ${J\subset\I}$, with ${k:=|J|\ge2}$, and let~$d_J$ be the degree of the greatest common divisor~$g$ of monomials in~$J$. Let $i$ be number of the variables that do not occur in~$g$, say $X_{\alpha_1},\ldots,X_{\alpha_i}$. Then~$J$ intersects exactly~$i$ faces of~$\F_{N,d}$, with ${0\le i\le N}$, since we are assuming~$J$ has all multiples of~$g$ that belong to~$\I$ (see remark~\ref{maxJ}). Now we know that the number of all monomials that are multiples of~$g$ is ${\tbinom{d-d_J+N}{N}}$. Since the multiples of $g$ that do not belong to any of the $i$ faces that intersect~$J$ can be written as $fgX_{\alpha_1}\cdots X_{\alpha_i}$, its number is ${\tbinom{d-d_J+N-i}{N}}$. Then exactly ${\tbinom{d-d_J+N}{N}-\tbinom{d-d_J+N-i}{N}}$ of its monomials are in~$\F_{N,d}$, and the remaining admit a greatest common divisor of degree ${d_J+i}$, and come from a subset~$J'\subseteq I'$, admitting a greatest common divisor of degree~${d_{J'}\ge d_J-N-1+i}$. Since ${N+1-i}$ faces of~$\F_{N,d}$ do not intersect~$J$, its greatest common divisor is multiple of the variables that are not present in those faces, and consequently ${d_J\ge N+1-i}$. Let
\[k':=\left|J'\right|=k-
    \left[\tbinom{d-d_J+N}{N}-\tbinom{d-d_J+N-i}{N}\right].\]
Observe that
\[\big[d'-(d_J-N-1+i)\big]n'+(d_J-N-1+i)-d'k'\ge
    \left(d'-d_{J'}\right)n'+d_{J'}-d'k'\mbox{,}\]
and the last expression is strictly positive, since~$I'$ satisfies inequality~(\ref{6.6}). We can see that
\begin{align*}
(d-d_J)n+d_J-dk&=\big[d'-(d_J-N-1+i)\big]n'+(d_J-N-1+i)-d'k'\\
    &\qquad +P(n',k',N,d,d_J,i)+Q(N,d,d_J,i)\mbox{,}
\end{align*}
where
\[P(n',k',N,d,d_J,i)=i\left(n'-k'\right)+
        (N+1-i)\left[\tbinom{d-d_J+N-i}{N}-k'+1\right]\]
and
\[Q(N,d,d_J,t)=(d-d_J)\left[\tbinom{d+N}{N}
        -\tbinom{d-1}{N}\right]-d\tbinom{d-d_J+N}{N}
        +(d-N-1+t)\tbinom{d-d_J+N-t}{N}.\]
Now~$P$ is clearly positive, since $\tbinom{d-d_J+N-i}{N}$ is the highest possible cardinality for~$J'$. If we can guarantee that $Q(N,d,d_J,i)$ is non\h{negative}, inequality~(\ref{6.6}) is strictly satisfied.

Suppose ${d_J=1}$. In this case, we have ${i=N}$, since ${d_J\ge N+1-i}$. Therefore we get ${Q(N,d,1,N)=\tfrac{(N-2)d+d-N}{d}\tbinom{d-1+N}{N}}$, and this expression is positive for ${N\ge3}$.

Now suppose ${d_J\ge2}$. Since the last term in~$Q$ vanishes for ${a>d-d_J}$, we shall consider the cases ${i\le d-d_J}$ and ${i>d-d_J}$ separately.

\bigskip
\noindent{\sc Case 1:} ${i\le d-d_J}$. Since
\[Q(N,d,d_J,t+1)-Q(N,d,d_J,t)=
    -\tfrac{(N-1)(d-N-2)+t(N+1)+(d_J-2)}{d-d_J-t}
    \tbinom{d-d_J+N-t-1}{N}\mbox{,}\]
and this is negative, we know that~$Q$ decreases as~$t$ gets higher. Therefore we should pay attention to the highest value of~$i$, that is ${i=\min(N,d-d_J)}$.

Let us start with the case ${d=N+2}$. In this case, ${i=N+2-d_J}$, and
\begin{align*}
Q(N,\,N+2,\,d_J,\,N+2-d_J)&=
    (N+2-d_J)\tbinom{2N+2}{N}-(N+2)\tbinom{2N+2-d_J}{N}\\
    &\qquad -N^2-2N+Nd_J+1.
\end{align*}
If we look at this as a function on~$d_J$, its second derivative is
\[-2\cdot\tfrac{N+2}{N!}\cdot\sum_{1\le s<t\le N}\
    \prod_{\substack{1\le r\le N\\t\neq r\neq s}}
    (N+2-d_J+r).\]
Since this is negative for ${d_J\in[2,N+1]}$, the function's minimum in this
interval is at one of its ends. For ${d_J=2}$, we have
\[Q(N,\,N+2,\,2,\,N)=\tfrac{3N^2-2N-4}{N+2}\tbinom{2N}{N}-N^2+1\mbox{,}\]
and this is positive for ${N\ge3}$. For ${d_J=N+1}$, we have
\begin{align*}
Q(N,\,N+2,\,N+1,\,1)&=\tbinom{2N+2}{N}-(N+2)\tbinom{N+1}{N}-N+1\\
    &\ge\tbinom{2N+2}{3}-(N+2)(N+1)-N+1\\
    &=\tfrac{1}{3}\left(4N^3+3N^2-10N-3\right)\mbox{,}
\end{align*}
and this is also positive for ${N\ge3}$.

\bigskip

Now, for general ${d\ge N+2}$, we shall see two subcases separately, namely\linebreak[4] ${N\le d-d_J}$ and ${N>d-d_J}$.

\bigskip
\noindent{\sc Subcase 1.1:} ${N\le d-d_J}$. Looking
at~$Q(N,d,d_J,N)$ as a function on~$d_J$ again, its second
derivative is
\[-\tfrac{2}{N!}\cdot\sum_{1\le s<t\le N}
    \left[d\cdot\sideset{}{_{t\neq r\neq s}}\prod_{r=1}^N
    (d-d_J+r)
    -(d-1)\cdot\sideset{}{_{t\neq r\neq s}}\prod_{r=1}^N
    (d-d_J-N+r)\right].\]
This is negative for ${d_J\in[2,d-N]}$, since ${d-d_J+r\ge d-d_J-N+r\ge1}$ for\linebreak[4] ${1\le r\le N}$. Therefore the function's minimum in this interval is at one of its ends. If ${d_J=2}$, we get
\begin{multline*}
Q(N,\,d,\,2,\,N)=\tfrac{1}{N!}(2d^2N+dN^2-2d^2-5dN-2N^2+2d+2N)
        \prod_{r=1}^{N-2}(d+r)\\
    -\tfrac{1}{N!}(N-1)\prod_{r=1}^{N}(d-N-1+r).
\end{multline*}
We will show that this expression is positive by induction on $N$. For ${N=3}$, we get
\[Q(3,\,d,\,2,\,3)=\tfrac{1}{3}d\left[d^2+6(d-4)+5\right]>0.\]
Now
\begin{align*}
Q(N+1,d,2,N+1)&=\tfrac{d+N-1}{N+1}Q(N,d,2,N)\\
    &\qquad+\tfrac{1}{(N+1)!}\left[2\left(d^2+dN-2d-2N\right)
        \prod_{r=1}^{N-1}(d+r)\right.\\
    &\left.\qquad -(d^2-2dN^2-d+2N^3-N^2+N)\prod_{r=1}^{N-1}(d-N+r)\right].
\end{align*}
Since
\begin{multline*}
2\left(d^2+dN-2d-2N\right)-\left(d^2-2dN^2-d+2N^3-N^2+N\right)=\\
    2(d-N)N^2+d(d-3)+2(d-1)N+(N-3)^2+3(N-3)>0,
\end{multline*}
and ${d+r>d-N+r}$ we get the result.

If ${d_J=d-N}$, we get
\[Q(N,\,d,\,d-N,\,N)=N\left[\tbinom{d+N}{N}-\tbinom{d-1}{N}\right]
    -d\tbinom{2N}{N}+d-1.\]
If we look at this expression now as a function on~$d$, its second derivative is
\[\tfrac{2}{(N-1)!}\cdot\sum_{1\le s<t\le N}
    \left[\sideset{}{_{t\neq r\neq s}}\prod_{r=1}^N(d+r)
    -\sideset{}{_{t\neq r\neq s}}\prod_{r=1}^N
    (d-N-1+r)\right]>0\mbox{,}\]
since ${d+r\ge d-N-1+r\ge1}$ for ${1\le r\le N}$. This means its first derivative increases with~$d$. When we evaluate this first derivative at ${d=N+2}$, we get
\begin{equation}\label{Brenner1}
\tfrac{1}{(N-1)!}\sum_{s=1}^N\left[\sideset{}{_{r\neq s}}
    \prod_{r=1}^N(N+2+r)-\sideset{}{_{r\neq s}}
    \prod_{r=1}^N(1+r)\right]-\tbinom{2N}{N}+1.
\end{equation}
Note that, except for the last term in the sum, the last factor of the first product is ${2(N+1)}$, and the last factor of the second is ${N+1}$. Splitting the last term in the sum and factoring out the last factor of both products, this becomes
\begin{multline*}
\tfrac{N+1}{(N-1)!}\sum_{s=1}^{N-1}\left[2\cdot
        \sideset{}{_{r\neq s}}\prod_{r=1}^{N-1}(N+2+r)
        -\sideset{}{_{r\neq s}}\prod_{r=1}^{N-1}(1+r)
        \right]\\
    +\tfrac{1}{(N-1)!}\left[\prod_{r=1}^{N-1}(N+2+r)
        -\prod_{r=1}^{N-1}(1+r)\right]-\tbinom{2N}{N}+1\mbox{,}
\end{multline*}
and this can be rearranged as
\begin{multline*}
\tfrac{N+1}{(N-1)!}\sum_{s=1}^{N-1}
    \left[\sideset{}{_{r\neq s}}\prod_{r=1}^{N-1}(N+2+r)
    -\sideset{}{_{r\neq s}}\prod_{r=1}^{N-1}(1+r)\right]\\
    +\tfrac{1}{(N-1)!}\left[\prod_{r=1}^{N-1}(N+2+r)
    -\prod_{r=1}^{N-1}(1+r)\right]\\
    +\tfrac{1}{N!}\left[\sum_{s=1}^{N-1}(N+1)N
    \sideset{}{_{r\neq s}}\prod_{r=1}^{N-1}(N+2+r)
    -\prod_{r=1}^{N}(N+r)\right]+1.
\end{multline*}
Now the first two terms are clearly positive and for the third one we can see that
\begin{multline*}
(N+1)N\sum_{s=1}^{N-1}\sideset{}{_{r\neq s}}
        \prod_{r=1}^{N-1}(N+2+r)-\prod_{r=1}^{N}(N+r)\\
    \ge(N+1)^2N\prod_{r=1}^{N-2}(N+2+r)-\prod_{r=1}^{N}(N+r)\\
    =(N+1)^2N\prod_{r=1}^{N-2}(N+2+r)-2N(2N-1)
        \prod_{r=1}^{N-2}(N+r)>0\mbox{,}
\end{multline*}
since ${(N+1)^2N>2N(2N-1)>1}$ for ${N\ge3}$, and ${N+2+r>N+r>1}$ for\linebreak[4] ${1\le r\le N-2}$. Therefore all the expressions are positive, and we see that~$Q$ increases with~$d$. Since we have already seen that~$Q$ is positive for the first value ${d=N+2}$, we have that~$Q$ is positive in this case.

\bigskip
\noindent{\sc Subcase 1.2:} ${N>d-d_J}$. Looking once more
at~$Q(N,d,d_J,N)$ as a function on~$d_J$, its second derivative is
\[-\tfrac{2d}{N!}\cdot\sum_{1\le s<t\le N}
    \sideset{}{_{t\neq r\neq s}}\prod_{r=1}^N
    (d-d_J+r).\]
This is negative for ${d_J\in[d-N,d-1]}$, since ${d-d_J+r\ge1}$ for ${1\le r\le N}$.\linebreak[4] Therefore the function's minimum in this interval is at one of its ends. If ${d_J=d-N}$, we are exactly in the same situation as before, so we already know that~$Q$ is non\h{negative}. If ${d_J=d-1}$, we get
\[Q(N,\,d,\,d-1,\,N)=\tbinom{d+N}{N}-\tbinom{d-1}{N}-d(N+1).\]
Therefore, if we apply lemma~\ref{Brenner2}, we get
\begin{align*}
Q(N,\,d,\,d-1,\,N)&\geq\tfrac{N+1}{(N-1)!}d^{N-1}-d(N+1)\\
    &\ge\tfrac{(N+1)^{N-1}}{(N-1)!}d-d(N+1)\\
    &>(N+1)d-d(N+1)=0.
\end{align*}

\bigskip
\noindent{\sc Case 2:} ${i>d-d_J}$. In this case, we get
${d_J>d-N}$, since ${i\le N}$. Let us start again with the case
${d=N+2}$. In this case, ${3\leq d_J\leq N+1}$, and
\[Q(N,\,N+2,\,d_J,\,i)=(N+2-d_J)\left[\tbinom{2N+2}{N}
    -\tbinom{N+1}{N}\right]-(N+2)\tbinom{2N+2-d_J}{N}.\]
If we look at this as a function on~$d_J$, its second derivative is
\[-2\cdot\tfrac{N+2}{N!}\cdot\sum_{1\le s<t\le N}\quad
    \sideset{}{_{t\neq r\neq s}}\prod_{1\le r\le N}
    (N+2-d_J+r).\]
Since this is negative for ${d_J\in[3,N+1]}$, the function's minimum in this interval is at one of its ends. For ${d_J=3}$, we have
\[Q(N,\,N+2,\,3,\,i)=\tfrac{7N^2-8N-8}{N+2}\tbinom{2N-1}{N}-N^2+1,\]
and this is positive for $N\ge3$. For ${d_J=N+1}$, we have
\begin{multline*}
Q(N,\,N+2,\,N+1,\,i)=\tbinom{2N+2}{N}-(N+3)(N+1)\\
    \ge\tbinom{2N+2}{3}-(N+3)(N+1)=\tfrac{1}{3}\left(4N^3+3N^2-10N-9\right)\mbox{,}
\end{multline*}
and this is also positive for ${N\ge3}$.

In general, for ${d\ge N+2}$, we get
\[Q(N,d,d_J,i)=(d-d_J)\left[\tbinom{d+N}{N}
    -\tbinom{d-1}{N}\right]-d\tbinom{d-d_J+N}{N}.\]
Looking at this expression as a function on~$d_J$, its second derivative is
\[-\tfrac{2d}{N!}\cdot\sum_{1\le s<t\le N}
    \sideset{}{_{t\neq r\neq s}}\prod_{r=1}^N(d-d_J+r)\le0.\]
Since this is negative for ${d_J\in[d-N,d-1]}$, the function's minimum in this interval is again at one of its ends.

For ${d_J=d-N}$, we have
\[Q(N,\,d,\,d-N,\,i)=N\left[\tbinom{d+N}{N}
    -\tbinom{d-1}{N}\right]-d\tbinom{2N}{N}.\]
If we look at this expression now as a function on~$d$, its second derivative is
\[\tfrac{2}{(N-1)!}\cdot\sum_{1\le s<t\le N}
    \left[\sideset{}{_{t\neq r\neq s}}\prod_{r=1}^N(d+r)
    -\sideset{}{_{t\neq r\neq s}}\prod_{r=1}^N
    (d-N-1+r)\right]>0\mbox{,}\]
since ${d+r\ge d-N-1+r\ge1}$ for ${1\le r\le N}$. This means its first derivative increases with~$d$. When we evaluate this first derivative at ${d=N+2}$, we get
\[\tfrac{1}{(N-1)!}\sum_{s=1}^N\left[\sideset{}{_{r\neq s}}
    \prod_{r=1}^N(N+2+r)-\sideset{}{_{r\neq s}}
    \prod_{r=1}^N(1+r)\right]-\tbinom{2N}{N}.\]
this is the same expression as~(\ref{Brenner1}) above, except for the final term ($+1$). Since this term is ignored in the argument that follows, we arrive to the same conclusions here, and guarantee that~$Q$ is positive.

For ${d_J=d-1}$, we have
\[Q(N,\,d,\,d-1,\,i)=\tbinom{d+N}{N}-\tbinom{d-1}{N}-d(N+1).\]
Again we can apply lemma~\ref{Brenner2}, and get
\begin{multline*}
Q(N,\,d,\,d-1,\,i)\geq\tfrac{N+1}{(N-1)!}d^{N-1}-d(N+1)\\
    >\tfrac{(N+1)^{N-1}}{(N-1)!}d-d(N+1)>(N+1)d-d(N+1)=0.
\end{multline*}

We have verified that $Q(N,d,d_J,i)$ is positive in all cases, so inequality~(\ref{6.6}) is strictly satisfied, and~$\I$ is a family of $n$ monomials whose associated syzygy bundle is stable.
\end{proof}

The main theorem of this work can now be stated. The fact that for the case ${N=2}$ problem~\ref{6.9} is solved (except for the case ${d=2}$ and ${n=5}$) \cite{CMMR10}, combined with the results in this section, will allow us to assert its main theorem. To get round that exception, we can see a particular case.

\begin{lem}\label{case326}
The syzygy bundle associated to the family of quadric monomials
${\I:=\big\{X_0^2,\,X_1^2,\,X_2^2,\,X_3^2,
    \,X_0X_1,\,X_2X_3\big\}}$
is stable.
\end{lem}
\begin{proof}
Note that the ideal generated by~$\I$ is primary, and that the relevant sets to verify inequality~(\ref{6.6}) have two elements and a linear greatest common divisor. Therefore
\[(d-d_J)n+d_J-dk=(2-1)\cdot6+1-2\cdot2=3>0\mbox{,}\]
so stability is guaranteed.
\end{proof}

\begin{te}\label{main3}
Let~$N$, $d$ and~$n$ be integers such that ${N\ge2}$, ${(N,d,n)\ne(2,2,5)}$, and ${N+1\le n\le\tbinom{d+N}{N}}$. Then there is a family of~$n$ monomials in $K\left[X_0,\ldots,X_N\right]$ of degree~$d$ such that their syzygy bundle is stable.

For ${(N,d,n)=(2,2,5)}$, there are~$5$ monomials of degree~$2$ in $K\left[X_0,X_1,X_2\right]$ such that their syzygy bundle is semistable.
\end{te}
\begin{proof}
Case ${N=2}$ was already stated in theorem 3.5 in \cite{CMMR10}.
For ${N\ge3}$, this can be done by induction on~$N$. For ${N=3}$,
lemma~\ref{part of a face and opposite vertex} gives us an answer
for
\[4\le n\le\tbinom{d+2}{2}+1\mbox{,}\]
except for the case where ${d=2}$ and ${n=6}$, for which we have lemma~\ref{case326}; proposition~\ref{faces} takes care of the cases
\[\tbinom{d+2}{2}+1<n\le\tbinom{d+3}{3}-\tbinom{d-1}{3}\mbox{;}\]
proposition~\ref{facesanddots} gives an answer for
\[\tbinom{d+3}{3}-\tbinom{d-1}{3}<n\le\tbinom{d+3}{3}
    -\tbinom{d-1}{3}+4\mbox{;}\]
finally, lemma~\ref{Brenner} takes care of all other cases.

Now if we suppose the answer is positive for some~$N$, lemma~\ref{part of a face and opposite vertex} provides a positive answer for the first cases of ${N+1}$, and proposition~\ref{faces}, proposition~\ref{facesanddots}, and lemma~\ref{Brenner} take care of the rest.
\end{proof}

\bibliographystyle{amsalpha}

\providecommand{\bysame}{\leavevmode\hbox to3em{\hrulefill}\thinspace}
\providecommand{\MR}{\relax\ifhmode\unskip\space\fi MR }
\providecommand{\MRhref}[2]{%
  \href{http://www.ams.org/mathscinet-getitem?mr=#1}{#2}
}
\providecommand{\href}[2]{#2}

\end{document}